\documentclass[12pt,a4paper]{amsart}
\usepackage{amssymb}
\usepackage{graphicx}  
\makeatletter
\renewcommand{\@begintheorem}[2]{
\rm \trivlist \item [\hskip \labelsep {\bf #2\ \ #1.}]
                                }
\makeatother

\makeatletter

\DeclareFontFamily{U}{cyr}{}
\DeclareFontShape{U}{cyr}{m}{n}{
  <5> wncyr5 <6> wncyr6 <7> wncyr7 <8> wncyr8 <9> wncyr9 <10->
wncyr10}{}
\DeclareMathAlphabet{\mathcyr}{U}{cyr}{m}{n}

\input cyracc.def
\input epsf

\newcommand\plotpicture{
\ 

\subsection{Figure}
Figure \ref{F:surface} shows the real locus of the surface given by $\det M_0=0$.
We depicted the part of affine chart with $x_3=1$ given by $-20 \leq x_0,x_1,x_2 \leq 20$. 
The figure also shows all points $g^n(x)$ in this region with $0\leq n <1000$
for $x = [0:1:0:1]$ (yellow), $x = [0:1:0:0]$ (red), and $x = [0:-1:0:1]$ (blue).

\begin{figure}[h]
\includegraphics[scale=0.5]{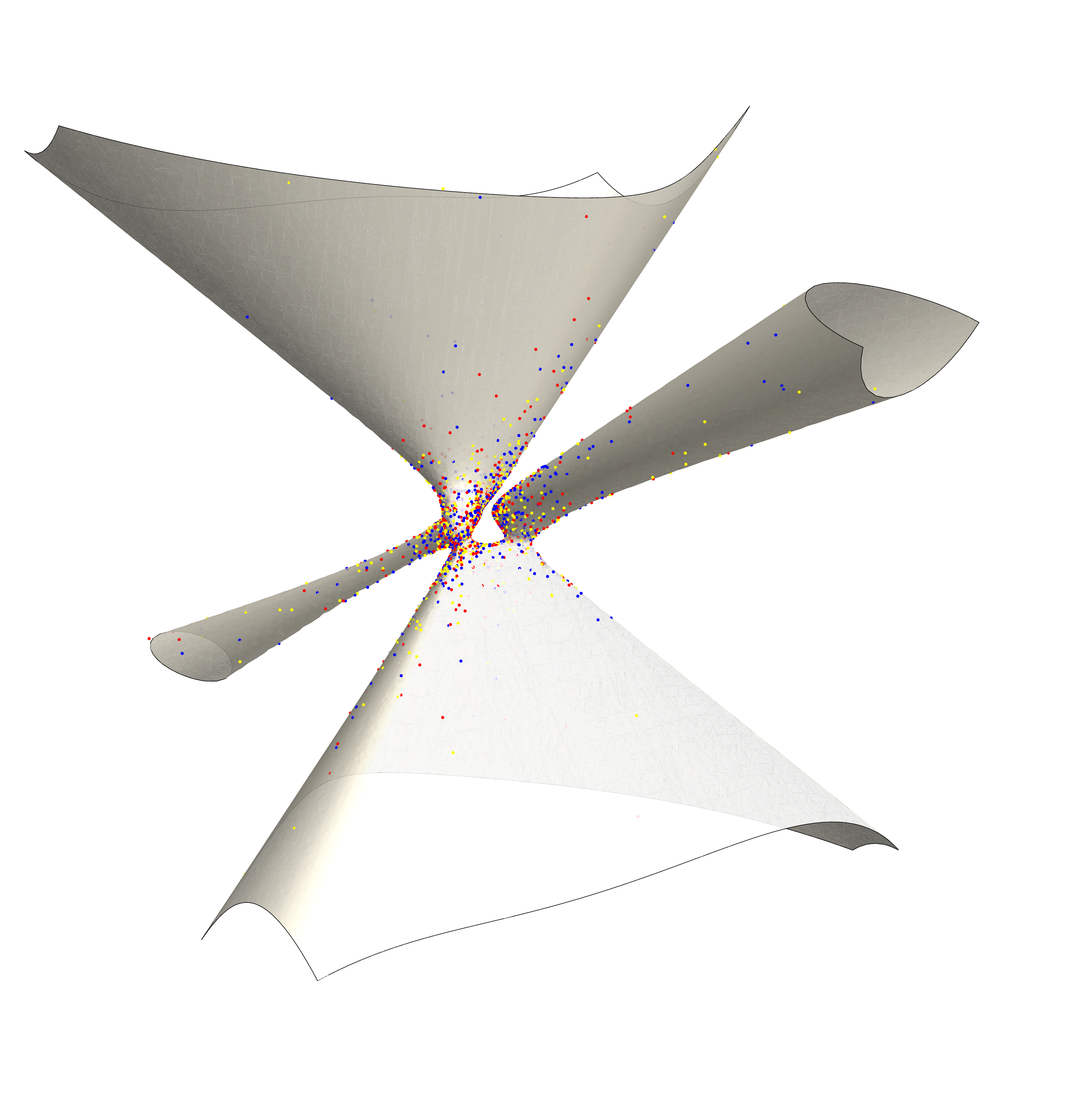}
\caption{The surface given by $\det M_0=0$ with points in three $g$-orbits.}
\label{F:surface}
\end{figure}
}

\newcommand\X{\mathfrak{X}}
\newcommand\p{\mathfrak{p}}
\newcommand\hH{\mathrm{H}}
\newcommand\Kbar{\overline{K}}
\newcommand\kbar{\overline{k}}
\newcommand\NS{\mathrm{NS}}
\newcommand\Aut{\mathrm{Aut}}
\newcommand\Spec{\mathop{\rm Spec}}
\newcommand\fS{\mathfrak{S}}
\newcommand\Tr{\mathop{\rm Tr}}
\newcommand\Gal{\mathop{\rm Gal}}

\newcommand{\ts}{\vspace{\baselineskip}\noindent{\bf Proof.}$\;\;$}
\newcommand{\ZZ}{{\Bbb Z}}
\newcommand{\QQ}{{\Bbb Q}}
\newcommand{\RR}{{\Bbb R}}
\newcommand{\CC}{{\Bbb C}}

\newcommand{\FF}{{\Bbb F}}

\newcommand{\PP}{{\Bbb P}}

\newcommand{\cA}{{\mathcal A}}

\newcommand{\cL}{{\mathcal L}}

\newcommand{\cO}{{\mathcal O}}

\newcommand{\inj}{\hookrightarrow}

\newcommand{\et}{{\text{\'et}}}

\newcommand{\im}{\mbox{im}}

\newcommand{\bes}{\begin{equation*}}
\newcommand{\ees}{\end{equation*}}

\title{The Cayley-Oguiso automorphism of positive entropy on a K3 surface}
\author{Dino Festi}
\author{Alice Garbagnati}
\author{Bert van Geemen}
\author{Ronald van Luijk}
\address{Matematisch Instituut, Universiteit Leiden,
Niels Bohrweg 1, 2333 Leiden, Nederland}
\address{Dipartimento di Matematica, Universit\`a di Milano,
Via Saldini 50, 20133 Milano, Italia}
\address{Dipartimento di Matematica, Universit\`a di Milano,
Via Saldini 50, 20133 Milano, Italia}
\address{Matematisch Instituut, Universiteit Leiden,
Niels Bohrweg 1, 2333 Leiden, Nederland}
\email{dinofesti@hotmail.it}
\email{alice.garbagnati@unimi.it}
\email{lambertus.vangeemen@unimi.it}
\email{rmluijk@gmail.com}
\begin{document}

\begin{abstract}
Recently Oguiso showed the existence of K3 surfaces that admit a 
fixed point free automorphism of positive entropy.
The K3 surfaces used by Oguiso have a particular rank two Picard lattice. 
We show, using results of Beauville, that these surfaces are therefore 
determinantal quartic surfaces. 
Long ago, Cayley constructed an automorphism of such determinantal surfaces.
We show that Cayley's automorphism coincides with Oguiso's free automorphism.
We also exhibit an explicit example of a determinantal quartic whose Picard lattice has exactly rank two and for which we thus have an explicit description of the automorphism.
\end{abstract}

\maketitle

Recently Keiji Oguiso showed that there exist projective K3 surfaces 
$S$ with a fixed point free automorphism $g$ of positive entropy, 
i.e.\ $g^*$ has at least one eigenvalue 
$\lambda$ of absolute value $|\lambda|>1$ on $H^2(S,\CC)$ (see \cite{O}).
He also described the Picard lattice of the general such surface explicitly and 
observed that these surfaces can be embedded into $\PP^3$ as quartic surfaces. 
There remained the problem of describing these quartic surfaces and their automorphism $g$ explicitly.

The aim of this paper is to provide a general method for 
constructing such quartic surfaces in $\PP^3$ and to describe an algorithm for finding the automorphism. 
Moreover, we will give an explicit example of such a surface $S$ and automorphism $g$.
To identify the quartic surfaces in Oguiso's construction, we observe that the Picard lattice required by Oguiso is exactly the Picard lattice of a general determinantal
quartic surface, that is, the quartic equation of the surface
is the determinant of a $4\times 4$ matrix of linear forms. 

While writing the paper, we realised that such automorphisms were already described 
by Prof.\ Cayley, President of the London Mathematical Society,
in his memoir on quartic surfaces, presented on February 10, 1870 (\cite{C}, $\S$69, p.47). 
In fact, Cayley observed that a determinantal K3 surface $S_0\subset \PP^3$
has three embeddings $S_i\subset\PP^3$ for $i=0,1,2$, 
each of which is again determinantal. 
The corresponding three matrices $M_i$,
which are closely related to each other,
provide natural (non-linear!) maps between these three quartic surfaces. 
A composition of these maps is an automorphism of $S_0$ and
we show that this automorphism is the one discovered by Oguiso. 

In the first section we recall Oguiso's description \cite{O} of K3 surfaces with 
a fixed point free automorphism $g$ of positive entropy. In section \ref{findg}
we give a method that in principle allows one to give an explicit description of the automorphism. In practice, even if the K3 surface $S$ is given as a determinantal surface in $\PP^3$, this method is hard to use, since one needs to know certain curves of high degree on $S$ that are not complete intersections.

In the second section, using results of Beauville, 
we give a characterisation of the K3 surfaces considered by Oguiso
as determinantal quartics.
Given the matrix $M_0(x)$ whose determinant is a defining polynomial for $S_0$, 
Cayley indicated a method to find the corresponding matrices $M_1(y)$ and $M_2(z)$ for $S_1$ and $S_2$ which we recall in Section \ref{cayleyA}.
We were not able to show that the determinants of $M_1$ and $M_2$ do not vanish identically in general. 
However, in the explicit example presented in Section \ref{explicit}, his method works and this allows us to give a convenient explicit description of the automorphism $g$ in that case.

After we put the first version of this paper on the arXiv, Igor Dolgachev informed 
us that the Cayley-Oguiso automorphism was (apparently independently) also discovered by F.\ Schur (\cite{S}, $\S$13, Satz I, p.30). The paper by Snyder and Sharp \cite{SS} presents the automorphism in a way similar to our Section \ref{cayleyA}. 
In a series of papers, \cite{R1},$\ldots$,\cite{R4}, T.G.\ Room studies the automorphism, especially in case the surface $S$ also contains a rational curve (in this case $\NS(S)$ has rank at least three). It is somewhat remarkable that none of these papers cites Cayley, but all refer to \cite{S} for the automorphism.

\

\noindent
{\bf Acknowledgements.} 
We are indebted to K.\ Oguiso and I.\ Dolgachev for helpful comments. 
We thank the referees for suggesting improvements to the first version of the paper.

\

\section{The general constructions}\label{general}

\subsection{The lattice $(N,b)$}\label{ns}
To describe the N\'eron Severi group of the K3 surfaces considered by Oguiso, we introduce a lattice $(N,b)$. 
One has $N\cong\ZZ^2$, but to describe $b$ and the isometries of $(N,b)$ it is convenient to define $N$ as the quotient ring $\ZZ[X]/(X^2-X-1)$ and let $\eta$ be the class of $X$: 
$$
N\,:=\,\ZZ[\eta],\qquad \eta^2\,=1\,+\,\eta~.
$$
The free $\ZZ$-module of rank two $N$ is isomorphic to the ring of integers of the number field
$\QQ(\eta)\stackrel{\cong}{\rightarrow}\QQ(\sqrt{5})$, $\eta\mapsto (1+\sqrt{5})/2$.
We denote the Galois conjugate of $x\in \QQ(\eta)$ by $x'$, 
so $(r+s\sqrt{5})'=r-s\sqrt{5}$ for $r,s\in \QQ$. One has 
$$
\eta'\,=\,1-\eta,\qquad (a\,+\,b\eta)'\,=\,a\,+\,b\eta'\,=\,(a+b)\,-\,b\eta\qquad(a,b\in\ZZ)~.
$$
The norm of $x\in N$ is defined as $Norm(x)=xx'$. We define a  
bilinear form
$$
b=b_N:\,N\,\times\,N\,\longrightarrow \,\ZZ,\qquad\mbox{by}\quad
b(x,y)\,=\,2(x'y+xy')~.
$$
So we get a lattice $(N,b)$:
$$
(N,b)\,
\stackrel{\cong}{\longrightarrow}\,\left(\ZZ^2,\;S_b=S_{b_N}\,:=\,\begin{pmatrix}4&2\\2&-4\end{pmatrix}\,\right)~,\qquad
a+b\eta\,\longmapsto \begin{pmatrix}a\\b\end{pmatrix}~.
$$
One easily verifies that 
$$
b(x,x)\,=\,4xx'\,=\,4(a^2+ab-b^2),\qquad(x=a+b\eta\,\in\,N)~.
$$

Due to the factor $4$, we have $b(x,y)\in 4\ZZ$, in particular, there are no $x\in N$ with $b(x,x)=\pm 2$. The equation $b(x,x)=0$ has only $x=0$ as solution, since if $a^2+ab-b^2=0$
with $b\neq 0$ then $(a/b)^2+(a/b)-1=0$, but this quadratic equation has no solution $a/b\in\QQ$.

As $\eta\eta'=-1$,  the map 
$$
N\,\longrightarrow\,N,\qquad x\,=\,a+b\eta\,\longmapsto \,
\eta^2 x=(a+b)\,+\,(a+2b)\eta
$$
is an isometry of the lattice $(N,b)$ with inverse $\eta^{-2}=2-\eta$. Composing this map with itself $n$ times gives an isometry which we denote simply by $\eta^{2n}$.

The isometries of $(N,b)$ are given by the maps $x\mapsto\pm \eta^{2k}x$ and 
$x\mapsto\pm \eta^{2k}x'$ with $k\in\ZZ$. 
To see this, we use that an isometry is given by a
$2\times 2$ matrix $M\in GL(2,\ZZ)$ on $\ZZ^2$ such that  ${}^tMS_bM=S_b$. 
Equivalently $M$ has integer coefficients, ${}^tMS_b=S_bM^{-1}$ and $\det M=\pm 1$. As $x\mapsto x'$ is an isometry with determinant $-1$, we need only consider the case $\det M=+1$.
One finds that $M$ must be a matrix with rows $(a,b),(b,a+b)$ 
and $1=\det M=a^2+ab-b^2$.
Thus the action of $M$ is the multiplication by $u=a+b\eta$ and $uu'=1$,
so $u$ is a unit in the ring of integers $\ZZ[\eta]$ of $\QQ(\sqrt{5})$.
This group of units is well-known to be $\{\pm \eta^m:m\in\ZZ\}$.
Thus $uu'=1$ implies that $x=\pm\eta^{2k}$ for some integer $k$.

\

With these definitions, Oguiso proved the following theorem, except for a refinement
which we prove here.

\subsection{Theorem (\cite{O}, Theorem 4.1)}\label{thmO}
There exist K3 surfaces $S$ with $\NS(S)\cong (N,b)$,
these form a dense subset of an 18-dimensional family of K3 surfaces.
The automorphism group $\Aut(S)$ of any such surface $S$ 
is isomorphic to $\ZZ$. Any generator of $\Aut(S)$
is a  fixed point free automorphism of positive entropy. 
Moreover, there is a generator $g$ of $\Aut(S)$ such that
$g^*=\eta^{6}$ on $\NS(S)\subset \hH^2(S,\CC)$ 
and $g^*=-1$ on the orthogonal complement $T(S)$ of $\NS(S)$ in $\hH^2(S,\CC)$.

\ts
In view of Theorem 4.1 in \cite{O}, all we need to prove is that 
$\Aut(S)\cong \ZZ$ and that one of its two generators acts as stated in the theorem.
An automorphism $\phi$ of $\Aut(S)$ is determined by its action $\phi^*$ on 
$\hH^2(S,\ZZ)$ (\cite{BHPV}, VIII, Corollary 11.2). 
The map $\phi^*$ preserves the sublattices $\NS(S)$, $T(S)$ of $\hH^2(S,\ZZ)$ and
$\phi^*$ is uniquely determined by its restriction to the sublattice (of finite index) $T(S)\oplus \NS(S)$ of $\hH^2(S,\ZZ)$.

Let $\phi^*_T$ be the restriction of $\phi^*$ to $T(S)$.
Then $\phi^*_T$ preserves the Hodge structure on $T(S)$ and thus it preserves the
two-dimensional subspace $T_t:=(\hH^{2,0}(S)\oplus \hH^{0,2}(S))\cap T(S)_\RR$ 
of $T(S)_\RR:=T(S)\otimes_\ZZ\RR$ as well as its orthogonal complement 
$T_a:=T_t^\perp$.
The intersection form is positive, respectively negative, definite on these spaces.
Thus the orthogonal groups $O(T_t)$, $O(T_a)$ are compact. 
As $\phi^*_T\in O(T(S))$, a discrete group, also lies in a compact group,
it lies in a finite set. In particular, $\phi^*_T$ has finite order. 

If $\phi^*_T\neq 1$, a suitable power of it will have a prime order $p$, let
$\tau$ be such a power of $\phi^*_T$.  The eigenvalues of  $\tau$ are thus $p$-th roots 
of unity and not all eigenvalues are equal to $1$.
If $\tau$ has eigenvalue $1$ on $T(S)$, then either the complexification of the sublattice $T(S)^\tau$ of $\tau$-invariants  or its orthogonal complement would contain 
$\hH^{2,0}(S)$. Thus either the orthogonal complement of $T(S)^\tau$
is contained in $\NS(S)$ or $T(S)^\tau$ itself is contained in $\NS(S)$.
Both cases contradict that $T(S)$ is $\NS(S)^\perp$.

Thus $\tau^{p-1}+\tau^{p-2}+\ldots+1=0$ on $T(S)$, hence also on the dual lattice
$T^*(S)\subset T(S)\otimes\QQ$  and therefore also on the discriminant group 
$T^*(S)/T(S)\cong \NS(S)^*/\NS(S)\cong (\ZZ/2\ZZ)^2\times(\ZZ/5\ZZ)$ 
(\cite{O}, Proposition 3.3).
In case $p\neq 2,5$ this leads to a contradiction: 
$\tau$ induces an automorphism of the subgroup $\ZZ/5\ZZ$ of $5$-torsion elements of the discriminant group.
As $\tau$ has order $p$ and $\Aut(\ZZ/5\ZZ)$ has order $4$,
$\tau$ must be the identity on $\ZZ/5\ZZ$.
%Hence $\tau(x)=x$ if  $x\in T^*(S)/T(S)$ has order five.
But then $0=(\tau^{p-1}+\ldots+1)x=px$ for all $x\in\ZZ/5\ZZ$, a contradiction.
In case $p=5$, one considers similarly the action of $\tau$ on the 2-torsion subgroup $(\ZZ/2\ZZ)^2$ of the discriminant group to get a contradiction.

Thus $\phi^*_T$ has order $2^k$ for some integer $k$. 
In case $k\geq 2$, there is thus an integer $m$ such that the restriction of 
$\psi:=(\phi^*)^m$ to $T(S)$ has order $4$.
As above, we can rule out that the restriction of $\psi$ to $T(S)$ 
has eigenvalues $\pm 1$, so $\psi^2+1=0$ on $T(S)$ and 
also on the discriminant group.
Thus $\psi$ acts as $x\mapsto \pm 2x$ on the subgroup $\ZZ/5\ZZ$ of the discriminant group.
The actions of $\psi$ on $T(S)$ and $\NS(S)$ are related through the equality $\psi_{T(S)^*/T(S)}=\phi_{\NS(S)^*/\NS(S)}$
under the natural isomorphism of discriminant groups 
$T(S)^*/T(S)\cong \NS(S)^*/\NS(S)\cong N^*/N$. 
Generators of $N^*/N$ are given in \cite{O}, Proposition 3.3(2), 
and we checked that 
the generators $x\mapsto x'$, $x\mapsto -x$ and $x\mapsto \eta^2 x$ of the isometry group all induce $x\mapsto -x$ on the subgroup $\ZZ/5\ZZ$ of 
the discriminant group. 
Therefore, $S$ cannot have an automorphism that has order four on $T(S)$. 

Thus we must have $\phi^*_T=\pm 1$ and, similar to the proof of Theorem 4.1 in \cite{O}, an automorphism of $S$ induces either 
$-1$ on $T(S)$ and $\eta^{12k+6}$ on $\NS(S)$ or 
$+1$ on $T(S)$ and $\eta^{12k}$ on $\NS(S)$ for some integer $k$. 
So $\phi$ acts as $g^{2k+1}$ or $g^{2k}$ on 
$T(S)\oplus \NS(S)$, where $g$ is the fixed point free automorphism found by Oguiso.
Therefore $\phi$ is a power of $g$. 
\qed

\subsection{Fibonacci numbers}\label{fib}
We will need to know the following values of $\eta^{2n}\in N$ explicitly:
$$
\eta^2\,=\,1+\eta,\qquad \eta^4\,=\,2+3\eta,\qquad \eta^6\,=\,5+8\eta~,
$$
as well as their inverses, with $\eta^{-2}=(\eta')^2=2-\eta$:
$$
\eta^{-2}\,=\,2-\eta,\qquad \eta^{-4}\,=\,5-3\eta,\qquad \eta^{-6}\,=\,13-8\eta~.
$$
The reader will notice the appearance of Fibonacci numbers 
(\cite{O}, Lemma 3.1): $\eta^{2n}=a_{2n-1}+a_{2n}\eta$ where $a_1=a_2=1$ and $a_{n+1}=a_n+a_{n-1}$ for $n\geq 1$ and 
$\eta^{-2n}=a_{2n+1}-a_{2n}\eta$.
In particular, one has
$$
\eta^{2n}\,+\,\eta^{-2n}\,=\,a_{2n+1}\,+\,a_{2n-1}\quad(\in\,\ZZ).
$$

\subsection{Topological Lefschetz numbers}\label{topLnum}
With $S$ and $g$ as in Theorem \ref{thmO}, 
the eigenvalues of $(g^n)^*$ on $\NS(S)$
are $\eta^{6n},\eta^{-6n}$, and $(g^n)^*$ acts as $(-1)^n$ on the 20-dimensional orthogonal complement $T(S)$ of $\NS(S)$ in $\hH^2(S,\ZZ)$.  
Notice that $g^n$ acts as the identity on the one-dimensional cohomology groups $\hH^i(S,\CC)$ for $i=0,4$, that $\hH^j(S,\CC)=0$ for $j=1,3$,
and that $\eta^{6n}+\eta^{-6n}=a_{6n+1}\,+\,a_{6n-1}$.
Thus the topological Lefschetz number of $g^n$ is
$$
T(S,g^n)\,:=\,\sum(-1)^itr(g^*|\hH^i(S,\CC))\,=\,
2\,+\,(-1)^n20\,+\,a_{6n+1}\,+\,a_{6n-1}~.
$$
In particular, the topological Lefschetz number of $g$ is $0$, a crucial step in the proof of Theorem \ref{thmO}. 
The topological Lefschetz number of $g^2$ is $22+a_{13}\,+\,a_{11}=344$, 
hence $g^2$ does have fixed points (cf. \cite{O}, Remark 4.3). 
Similarly, if $|n|>1$ then $(g^n)^*$ has fixed points, but it has no fixed curves 
since $(g^n)^*$ has no eigenvalue $1$ on $\NS(S)$. 
In Section \ref{explicit} we will present an example where $g^2$ 
has exactly $344$ fixed points (see Proposition~\ref{periodic}).

\subsection{Ample divisors on $S$}\label{quartics}
Let $S$ be a K3 surface with Picard lattice $\NS(S)\cong (N,b)$ 
as in section \ref{ns}.
We will fix the identification $\NS(S)\cong N$ in such a way that if $D$ is an 
ample divisor class, so $D^2>0$, then $D=a+b\eta$ with $a>0$.
As there are no elements with $b(x,x)=-2$ in $N$, any $x=a+b\eta \in N$ with 
$b(x,x)>0$ and $a>0$ is the class of an ample divisor on $S$ 
(\cite{BHPV}, VIII, Corollary (3.9)):
$$
\cA(S)\,=\,
\big\{\,x=a+b\eta\,\in \NS(S):\;a>0,\quad b(x,x)\,=\,4(a^2+ab-b^2)>0\,\big\}~.
$$
The isometry $\eta^2$ of $\NS(S)$ maps $\cA$ onto itself: $\eta^2\cA=\cA$.
This is easily seen by observing that an isometry of $N$ extends $\RR$-linearly 
to an isometry of $N_\RR=\RR^2$ which maps the set $Q$
defined by $a^2+ab-b^2=0$ into itself.
This set consists of  two lines and the four connected components of $N_\RR-Q$ are thus permuted by an isometry.
The isometry $\eta^2$ maps $1\in\cA$ to $\eta^2\in\cA$, so it fixes the connected component containing $\cA$, hence $\eta^2\cA=\cA$.

\subsection{Ample divisors on $S$ are very ample}\label{ampvamp}
We recall the results of Saint-Donat which imply that any ample divisor on $S$ 
is already very ample (cf.\ \cite{O}, Remark 4.2). 
Let $L$ be the line bundle on $S$ defined by an ample divisor class $D$. 
As the canonical bundle of $S$ is trivial, 
Kodaira vanishing implies that $h^i(L)=\dim \hH^i(S,L)=0$ for $i>0$.
The Riemann-Roch theorem then asserts that $h^0(L)=2+D^2/2$. 
As an ample divisor $D$ on $S$ has $D^2>0$, 
we get $h^0(D)>0$ and thus we may assume that $D$ is effective.
As $\Delta^2\neq -2$ for any divisor $\Delta$ on $S$, the linear system $|D|$ has no fixed components (\cite{SD}, $\S$ 2.7.1, 2.7.2). 
By \cite{SD}, $\S$ 4.1, the map $\phi_L$ defined by the global sections of $L$ 
is then either of degree two or it is birational onto its image. 
In the first case, Theorem 5.2 of \cite{SD} implies that $S$ has 
a divisor $\Delta$ with $\Delta^2\in \{0,2\}$, which is not the case. 
So $\phi_L$ is birational onto its image and by Theorem 6.iii of \cite{SD}, 
the image of $\phi_L$ is isomorphic to the image of the map $\theta_L$
introduced in \cite{SD}, $\S$ 4.2. 
The map $\theta_L$ is the contraction of all smooth rational curves $\Delta$ in $S$ with $D\cdot\Delta=0$. 
As $\Delta^2=-2$ but there are no divisor classes with $\Delta^2=-2$ in $S$, 
the map $\theta_L$ is the identity on $S$ and $\phi_L:S\rightarrow \phi_L(S)$ is an isomorphism, 
i.e.\ $L$ is very ample.

\subsection{Effective divisors and irreducible curves}\label{effirr}
Let $D\in \NS(S)$ be the class of an irreducible curve $C$, then by adjunction $D^2=2p_a(C)-2\geq -2$, and thus actually $D^2>0$. As also $D\cdot H>0$ for any ample divisor $H$, we conclude that $D=a+b\eta$ with $a>0$ and therefore any curve in $S$ is an ample divisor. Taking linear combinations with positive coefficients of classes of curves, we conclude that any effective divisor on $S$ is an ample divisor.

%Let $D_0\in \NS(S)$ be the class with $a=1$, $b=0$, so $D_0=1$.
%The $D_0\in\cA$, and $D_0^2=4$, hence $h^0(D_0)=2+D_0^2/2=4$ 
%and curves in $|D_0|$ have genus $3$. 
%We claim that any divisor $D\in |D_0|$ is irreducible and reduced. 
%If not, we would have $D=D'+D''$, a sum of effective divisors, and thus %$D',D''\in\cA$. Then $D'=a'+b'\eta$, $D''=a''+b''\eta$ with $a',a''\in\ZZ$ and %$a',a''>0$, which contradicts that $D'+D''=1$.

\subsection{Quartic surfaces}\label{quas}
For any integer $n$ we define a very ample divisor class
$$
D_n\,:=\,\eta^{2n}\,\in N,\qquad D_n^2\,=\,D_0^2\,=\,4\qquad (n\in\ZZ)~.
$$
A basis of the global sections of the line bundle on $S$ defined by $D_n$
defines a projective embedding, denoted by $\phi_n$,
of $S$ in $\PP^3$ as a quartic surface $S_n\subset\PP^3$:
$$
\phi_n\,:=\,\phi_{D_n}:\,S\,\stackrel{\cong}{\longrightarrow}\,S_n\,\subset\,\PP^3~.
$$
%The divisors in $|D_n|$ are all irreducible and reduced, in fact if 
%$D_n=D_n'+D_n''$ with effective divisors
%$D_n',D_n''$, then $D_n',D_n''\in\cA$.
%Applying the isometry  $\eta^{-2n}$ to the classes of these divisors in 
%$N=\NS(S)$ by $\eta^{-2n}$, we obtain
%$D_0=\eta^{-2n}D_n'+\eta^{-2n}D''_n$ and, as $\eta^{2}\cA=\cA$, 
%we have $\eta^{-2n}D_n',\eta^{-2n}D''_n\in\cA$, but this impossible 
%as we showed in section \ref{effirr}.

\subsection{The automorphism $g$}\label{autg}
The remarkable fact that there is an automorphism of $S$ with $g^*=\eta^6$
implies that the quartic surfaces $S_0$ and $S_3$ are the same, 
after choosing suitable coordinates on the $\PP^3$'s.
 
In fact, let $s_0,\ldots,s_3$ be a basis of $\hH^0(S,D_0)$. 
As $g^*D_0=D_3$,
$\hH^0(S,D_3)$ has basis $t_i:=g^*s_i$, where $i=0,\ldots,3$. 
With a slight abuse of notation, we then get for all $x\in S$:
$$
\phi_{3}(x)\,=\,(t_0(x):\ldots:t_3(x))\,=\,(s_0(g(x)):\ldots:s_3(g(x)))\,=\,\phi_0(g(x))~.
$$
Thus, with these bases, $S_3=S_0\subset\PP^3$. 
Moreover $\phi_3=\phi_0\circ g$ implies that
$$
g\,=\,\phi_0^{-1}\circ\phi_3:S\,\longrightarrow\, S~.
$$

\subsection{How to find $g$}\label{findg}
To give a more concrete description of $g$, 
we explain how, in principle, one can describe $\phi_3$ in terms of $\phi_0$.
For this we need to find $\hH^0(S,D_3)$, given the surface $S_0\subset\PP^3$.
The zero locus of a global section $t$ of $D_3$ is mapped to a curve in $S_0$.
This curve is not the (complete) intersection of $S_0$ with another surface (of degree $d$) in $\PP^3$, since such an intersection has class $dD_0=d$, whereas $D_3=\eta^6=5+8\eta$. 

The intersection of such a surface of degree $d$ is thus the sum of two effective divisors with classes
$D_3$, $D$ respectively and $dD_0=D_3+D$. 
As effective classes and ample classes coincide on $S$,
the smallest possible degree $d$ is the smallest positive integer such that $dD_0-D_3\in\cA$, that is, $d-5>0$ and $(d-5)^2+(d-5)(-8)-(-8)^2>0$, which is $d=18$ (and then $D=18D_0-D_3=D_{-3}$).

%Notice that $\eta^{-6}=13-8\eta$ and thus $\eta^6+\eta^{-6}=18$, so
%$D_3+D_{-3}=18D_0$.
Let $C_i:=(t_i=0)$ be the zero divisors of a basis $t_i$, 
$i=0,\ldots,3$, of $\hH^0(S,D_3)$ and similarly, let
the  $t'_j$ be a basis of $\hH^0(S,D_{-3})$ with zero divisor $C'_j:=(t'_j=0)$. 
Then the divisor $C_i+C'_j$ has class $D_3+D_{-3}=18D_0$ and it
is the zero locus of the section $t_it_j'$ in $\hH^0(S,18D_0)$.

Consider the exact sequence of sheaves on $\PP^3$:
$$
0\,\longrightarrow\,\cO_{\PP^3}(d-4)\,\longrightarrow\,\cO_{\PP^3}(d)\,
\longrightarrow\, i_*\cO_{S_0}(d)\,\longrightarrow\,0~,
$$
where the first non-trivial map is multiplication by the equation of $S_0$ and 
where $i:S_0\hookrightarrow \PP^3$ is the inclusion map. 
As the first cohomology group $H^1(\PP^3,\cL)$ of any invertible sheaf $\cL$ on $\PP^3$ is zero, 
and $\phi_0^*\cO_{\PP^3}(d)=dD_0$,
we get a surjection 
$$
\hH^0(\PP^3,\,\cO_{\PP^3}(d))\,
\stackrel{\phi_0^*}{\longrightarrow}
\,\hH^0(S,dD_0)\,\longrightarrow\,0~.
$$
Therefore, for any $d$, a section in
$\hH^0(S,dD_0)$ is the restriction of a homogeneous polynomial of 
degree $d$ on $\PP^3$.

In particular, there are homogeneous polynomials $R_{ij}$ of degree $18$ in 
$x_0,\ldots,x_3$, such that 
$$
{\phi_0^*R_{ij}}\,=\,t_it_j'\quad\in \hH^0(S,18D_0)~,\qquad (i,j\in\{0,\ldots,3\})~.
$$
Considering the zero loci of these sections we get:
$$
(R_{ij}=0)\,\cap\, S_0\;=\,\phi_0(C_i)\,+\,\phi_0(C'_j)~.
$$
The curves $\phi_0(C_i)$, $\phi_0(C'_j)$ in $\PP^3$ both
have degree $D_0\cdot C_i=36=D_0\cdot C'_j$, 
consistent with $18\cdot 4=72=36+36$.

The map $\phi_3$ is defined by the global sections $t_0,\ldots,t_3$ of $D_3$.
Since $D_{-3}$ is very ample, for each $x\in S$ there is an index $j$ such that $t'_j(x)\neq 0$. So (with slight abuse of notation):
{\renewcommand{\arraystretch}{1.3}
$$
\begin{array}{crrrl}
\phi_3:\,S\,\longrightarrow\,S_3\,\subset\,\PP^3,\qquad&
p&\longmapsto &&(t_0(p):\;\ldots\;:t_3(p))\\
&&&=&(t_0(p)t'_j(p):\;\ldots\;:t_3(p)t'_j(p))\\
&&&=&(R_{0j}(\phi_0(p)):\ldots:R_{3j}(\phi_0(p)))~.
\end{array}
$$
} %end stretch
On the open subset of $S$ where $t'_j\neq 0$, we thus have: 
$\phi_3=R_j\circ\phi_0$, where $R_j:\PP^3\rightarrow\PP^3$
is the rational map given by the polynomials $R_{0j},\ldots,R_{3j}$. 
Hence on this open subset we get 
$$
g\,=\,\phi_0^{-1} \circ \phi_3\,=\,\phi_0^{-1}\circ R_j\circ\phi_0~,
$$
that is, if we identify $S$ with $S_0$, then $g$ is just the rational map $R_j$, 
for any $j$, and these maps glue to give an isomorphism $S_0\rightarrow S_0$,
which `is' $g$. 

To find $\phi_3$, given $\phi_0$, we thus need to find these polynomials $R_{ij}$ 
on $\PP^3$. 
In practice, even if one is given that $\NS(S)\cong (N,b)$ this seems quite difficult. 
However, in Section \ref{cayleyA}, 
we construct explicit polynomials of degree $27$ which induce the map $S_0\rightarrow S_0$,  
corresponding to $g$, using a general method due to Cayley. 
Using these polynomials, we were able to find the degree $18$ polynomials in the specific example in Section \ref{explicit},
see Section \ref{deg18}.

\

\section{Determinantal quartic surfaces}\label{detquars}

\subsection{Determinantal quartics}\label{detqua}
We now show that a result of Beauville provides an explicit description of the
K3 surfaces we are interested in:
the K3 surfaces $S$ with N\'eron Severi group isomorphic to $(N,b)$ are 
exactly the quartic determinantal surfaces with Picard number two. 

More precisely, the quartic surfaces $S_n:=\phi_n(S)$ from Section \ref{quas}
are all determinantal. In Corollary \ref{class} we show how the matrix $M_n$ which defines $S_n$ also provides explicitly the map 
$\phi_{n+1}\phi_n^{-1}:S_n\rightarrow S_{n+1}$, once 
suitable bases of global sections are chosen. 
This is actually part of the results of Cayley in \cite{C}.

\subsection{Proposition} \label{propB}
Let $S$ be a K3 surface with N\'eron Severi group 
$\NS(S)\cong (N,b)$ as in Section \ref{ns}. 
Then, for any $n\in\ZZ$, the quartic surface $S_n:=\phi_n(S)$ is determinantal.
So there is a $4\times 4$ matrix $M_n(x)$, whose coefficients are linear forms in 4 variables $x_0,\ldots,x_3$, such that $\det M_n(x)=0$ is an equation for $S_n$.

Conversely, a general determinantal quartic surface $S$ 
%has rank$(\NS(S))\geq 2$.
%A general  determinantal quartic surface $S$ has 
has $\NS(S)\cong (N,b)$
and thus it admits a fixed point free automorphism of positive entropy.

\ts
The proposition is an easy consequence of
\cite{B}, Proposition 6.2, where Beauville proved that a smooth quartic surface $X$ is determinantal if and only if there is a curve $C\subset X$ of degree $6$ and genus $3$. 
See also \cite{topics}, Section 4.2.5.

Given a K3 surface  $S$ with  $\NS(S)\cong (N,b)$, there are smooth genus three curves $C_n$ on $S$ with class the very ample divisor $D_n=\eta^{2n}$, for all 
$n\in\ZZ$ (cf.\ Section \ref{quartics}).
As multiplication by $\eta^{-2n}$ is an isometry of the lattice $(N,b)$, 
and as one easily computes $C_0\cdot C_1=6$, 
we then get that $C_n\cdot C_{n+1}=6$. 
As $\phi_n(C_n)$ is a plane section of $S_n$, 
the curve $\phi_n(C_{n+1})$ is a genus $3$ curve of degree $6$ in $S_n$. 
Hence $S$ is determinantal.

For the converse, let $S$ be a general determinantal quartic surface in $\PP^3$.
Then $S$ is smooth (\cite{B} (1.10)).
Let $H$ be the hyperplane class of $S$, so $H^2=4$. 
Let $C\subset S$ be a degree $6$ and genus $3$ curve as in  
\cite{B} Proposition 6.2. 
Then $H\cdot C=6$ and the adjunction formula implies that $C^2=4$.
Thus the intersection form on the sublattice $\ZZ H\oplus \ZZ C$ of $\NS(S)$ 
is given by the matrix
$$
\begin{pmatrix} H^2&H\cdot C\\ H\cdot C&C^2\end{pmatrix}
\,=\,
\begin{pmatrix} 4&6\\6&4\end{pmatrix}~.
$$
This sublattice is isometric to $(N,b)$ since the $\ZZ$-basis of $N$ given by $D_0=(1,0)$ and $D_1=(1,1)$ gives this intersection matrix. 
Thus $\NS(S)$ of a determinantal quartic K3 surface contains $(N,b)$ as a sublattice and therefore the rank of $\NS(S)$ is at least two.
%To show that $\NS(S)$ has rank two, we can use a computation of 
%Cayley (\cite{C}, $\S$66, p.45)
%which shows that determinantal quartics have $18$ moduli, 
%and thus the rank of the N\'eron Severi group is $2$ 
%for a general determinantal quartic. 
%Alternatively, 
In the next section we provide an example of a (smooth) determinantal quartic with $\mbox{rank}\,\NS(S)=2$, 
thus the same is true for the general determinantal quartic.

For such a quartic we thus have $N\subset \NS(S)$, of finite index. 
As $|\det(b)|=20$, the index can only be $1$ or $2$. If the index is two then $D:=(aH+bC)/2\in \NS(S)$ with $(a,b)=(1,0)$ or $(0,1)$ or $(1,1)$, but $D^2$ is odd in all these cases, so this is impossible. 
Hence $\NS(S)=(N,b)$ for a general  determinantal quartic surface.
The existence of a fixed point free automorphism of positive entropy now follows from Oguiso's results in \cite{O}.
\qed

\subsection{Generators of $\NS(S)$}\label{gens}
The proposition implies in particular that a general smooth determinantal surface $S\subset\PP^3$ has N\'eron Severi group of rank two. 
One would thus like to see a curve on $S$
which is not a complete intersection, that is, 
whose class is not an integer multiple of the hyperplane class $H$ of $S\subset\PP^3$. 
As explained in \cite{B} (see also \cite{topics}, Example 4.2.4), 
such curves, of genus 3 and degree 6, can be found as follows. 

The matrix of linear forms $M$, whose determinant defines $S$, also gives
a sheaf homomorphism $\cO(-1)^{\oplus 4}\rightarrow \cO^{\oplus 4}$
on $\PP^3$. 
The cokernel is $i_*\cL$ for an invertible sheaf $\cL$ on $S$, 
where $i:S\hookrightarrow \PP^3$ 
is the inclusion (\cite{B}, Corollary 1.8). 
$$
0\,\longrightarrow\,\cO(-1)^{\oplus 4}\,\stackrel{M}{\longrightarrow}\,
\cO^{\oplus 4}\,\longrightarrow\,i_*\cL\,\longrightarrow\,0~.
$$
So $M$ defines a line bundle on $S$ with sheaf of sections $\cL$. 
We will denote this line bundle by $\cL$ as well.
As $\hH^i(\PP^3,\cO(-1))=0$ for all $i$, we obtain an isomorphism
$$
\CC^4\,=\,\hH^0(\PP^3,\cO^{\oplus 4})\,\stackrel{\cong}{\longrightarrow}\,
\hH^0(S,\cL)~.
$$
In Proposition \ref{line} we will show that $\cL$ has sections whose zero locus has degree $6$ and genus $3$.
%We denote by $s_i\in\Gamma(S,\cL)$ 
%the global section defined by the $i$-th basis vector of $\CC^4$. 
\

\subsection{The cofactor matrix} \label{cofac}
We recall some well-known linear algebra.
For an $n\times n$ matrix $M=(m_{ij})$,
with coefficient $m_{ij}$ in the $i$-th row and $j$-th column, let $M_{ij}$ be the 
$(n-1)\times(n-1)$ matrix obtained from $M$ by deleting the $i$-th row and $j$-th column. The cofactor matrix of $M$ is the $n\times n$ matrix
$$
P\,:=\,(p_{ij})\qquad\mbox{with}\quad
p_{ij}\,:=\,(-1)^{i+j}\det(M_{ji})~.
$$
Let $I$ be the $n\times n$ identity matrix, then we have the following matrix identities:
$$
MP\,=\,PM\,=\,\det(M)I~.
$$

\subsection{Proposition}\label{line}(\cite{B}, (6.7))
Let $S$ be a smooth quartic surface defined by $\det M=0$. 
Let $\cL$ be the line bundle on $S$ defined in Section \ref{gens},
let $s_j\in \hH^0(S,\cL)$ be the global section of $\cL$ 
which is the image of the $j$-th basis vector of $\CC^4$ and 
let $C_j$ be the zero locus of $s_j$. 

Then $C_j\subset S$ is the divisor defined by the vanishing of the four
coefficients $p_{ij}$, $i=1,\ldots,4$ of the cofactor matrix $P$ of $M$. 
Moreover, the effective divisors $C_j$ have degree 6 and genus 3.

\ts As $S$ is smooth, for any $x\in S$ at least one of the partial derivatives
$(\partial/\partial x_i \det M)(x)\neq 0$. 
Using the expansion of the determinant of a matrix $M=(m_{ij})$,
whose coefficients are variables $m_{ij}$, according to the $i$-th row one finds
that $\partial/\partial m_{ij} \det M=(-1)^{i+j}\det M_{ij}$. 
As each $m_{ij}(x)$ is a  function of $x_0,\ldots,x_3$ one finds,
using the chain rule, that at least one 
$3\times 3$ minor $\det M_{ij}(x)$ is non-zero for $x\in S$.
%Thus $M(x)$ has rank $3$ for any $x\in S$.

Notice that $x\in C_j$ if and only if $e_j(x)\in \im M(x)$ where $e_j$ 
is the global section of the trivial bundle $\cO^{\oplus 4}$ 
defined by the $j$-th basis vector. 
%For $x\in S$ one has $\mbox{rank}\,M(x)=3$ since $S$ is smooth
%(to see this, use the expansion of the determinant of the matrix $X=(x_{ij})$
%whose coefficients are variables $x_{ij}$ according to the $i$-th row to obtain
%$\partial/\partial x_{ij} \det X=(-1)^{i+j}\det X_{ij}$
%and use the chain rule for differentiation).
Let $P(x)$ be the cofactor matrix of $M(x)$, then $P(x)M(x)=0$, 
which implies that $\im M(x)\subset \ker P(x)$. 
As at least one $3\times 3$ minor of $M(x)$ is non-zero, we have $P(x)\neq 0$ and we conclude that $\dim \ker P(x)=3$
and so $\im M(x)= \ker P(x)$. 
Thus $e_j(x)\in \im M(x)$ is equivalent to $P(x)e_j(x)=0$ 
which is equivalent to the vanishing of $p_{ij}(x)$
%=(-1)^{i+1}\det(M_{1i}(x))$, 
for $i=1,\ldots,4$.
%Up to sign, these are the maximal minors as in the proposition, 
%hence we get the first part of the proposition.

The degree and genus of $C=C_j$ are given in 
\cite{B} Prop.\ 6.2, \cite{topics} Thm.\ 4.12.14.
The genus is actually easy to compute in this case:
as $\cL\cong \cO_S(C)$ and all effective divisors on $S$ are ample, 
Kodaira vanishing and Riemann-Roch on $S$ imply that
$4=\dim \hH^0(S,\cO_S(C))=\chi(\cO_S(C))=p_a(C)+1$.
\qed

\

\subsection{The transposed matrix}\label{transsec}
Given a determinantal surface $S$ with equation $\det M=0$, 
one obviously also has the (same) equation $\det{}^tM=0$. 
However, the invertible sheaf $\cL'$ on $S$ defined by the cokernel of ${}^tM$
is not isomorphic to $\cL$, but to $\cO_S(3)\otimes\cL^{-1}$ (\cite{B}, (6.3),
\cite{topics}, (4.19)). 

%The following proposition, applied to ${}^tM$ instead of $M$, 
%thus shows that the map defined by the global sections of 
%$\cO_S(3)\otimes\cL^{-1}$ is defined by the rows of the cofactor matrix 
%of ${}^tM$, that is by the columns of the cofactor matrix of $M$.

\subsection{Proposition}\label{maps} 
Let $S$ be a determinantal surface defined by $\det M=0$,
let $\cL$ be the line bundle on $S$ defined in Section \ref{gens} 
and let
$s_1,\ldots,s_4$ be the basis of $\hH^0(S,\cL)$ 
defined in Proposition \ref{line}.
The rational map
$$
\phi_\cL:\,S\,\longrightarrow\,\PP^3~,\qquad 
x\,\longmapsto\,(s_{1}(x):\ldots:s_{4}(x))~,
$$
coincides with the rational map given by any of the rows of the cofactor matrix $P$ of $M$:
$$
S\,\longrightarrow\,\PP^3~,\qquad 
x\,\longmapsto\,(p_{i1}(x):\ldots:p_{i4}(x))~,
$$
for any $i$ in $\{1,\ldots,4\}$.

\ts 
In Proposition \ref{line} we showed that the coefficients $p_{ij}$, $i=1,\ldots,4$, 
of the cofactor matrix $P$ of $M$ define the zero locus $C_j$ of 
the section $s_j$ of $\cL$. 
As ${}^tP$ is the cofactor matrix of ${}^tM$, the coefficients $p_{ij}$, $j=1,\ldots, 4$, of the cofactor matrix $P$ similarly define the zero locus $C'_i$ of a section $t_i$ of $\cO_S(3)\otimes\cL^{-1}$.
In particular, the coefficient $p_{ij}$ of $P$ is zero on both $C_j$ and $C'_i$. 
%The zero locus of $p_{11}$ in $S$ is in fact exactly $C\cup C'$,
%since 
More precisely, we have the following identity of divisors on $S$:
$$
S\,\cap\,(p_{ij}\,=\,0)\;=\;C_j\,+\,C'_i~.
$$
In fact, $p_{ij}=0$ is a cubic surface in $\PP^3$ and thus the left hand side is a divisor with class
$3H$ in $\NS(S)$, where $H$ is the hyperplane class of $S$. The classes of the line bundles $\cL$, $\cO_S(3)\otimes\cL^{-1}$ are the classes of the divisors of the zero loci of their divisors $C_j$, $C'_i$ respectively.
But the class of the line bundle $\cO_S(3)\otimes\cL^{-1}$ is also $3H-C_j$, hence $C_j+C'_i=3H$, which proves the identity.
Now we define sections $t_i$ of $\cO_S(3)\otimes\cL^{-1}$, 
with zero locus $C'_i$, by $p_{i1}=t_is_1$.
As $P(x)$ has rank one for $x\in S$, we get $p_{ij}(x)p_{11}(x)=p_{i1}(x)p_{1j}(x)$ hence $p_{ij}=t_is_j$ for any $i,j$ and the proposition follows. 
\qed

\subsection{Explicitly moving from $S_n$ to $S_{n+1}$} 
Now we return to the quartic surfaces $S_n=\phi_n(S)$, where $S$ 
is a K3 surface with $\NS(S)\cong (N,b)$ as in Section \ref{general}.

Each $S_n$ is a determinantal surface, with equation $\det M_n=0$, by Proposition \ref{propB}.
The following corollary identifies the line bundle $\cL$ (up to replacing $M_n$ by ${}^tM_n$), it is the line bundle defined by the divisor class $D_{n+1}$ (or $D_{n-1}$).  
In particular, it basically solves the problem of moving from $S_n$ to $S_{n+1}$
as the following corollary shows.

%We need to know the classes of the two divisors $C$, $C'$ cut out by a 
%$3\times 3$ minor $p_{n,11}$, and more generally of $p_{n,ij}$, of $M_n$.
%This will allow us to apply the method we outlined in Section \ref{move} 
%to move from $S_n$ to $S_{n+1}$.

\subsection{Corollary} \label{class}
Let $S$ be a K3 surface with N\'eron Severi group 
$\NS(S)\cong (N,b)$ as in Section \ref{general}. 
Let $S_n=\phi_n(S)\subset \PP^3$ be the smooth determinantal surface 
defined by $\det M_n=0$ and 
let $P_n=(p_{n,ij})$ be the cofactor matrix of $M_n$.
Let $\cL$ be the line bundle on $S_n$ defined by $M_n$ as in Section \ref{gens}.

Then the line bundle $\phi_n^*\cL$ on $S$ has class $D_{n-1}$ 
(and $\phi_n^*(\cO(3)\otimes \cL^{-1})$ has class $D_{n+1}$)
or it has class $D_{n+1}$ 
(and then $\phi_n^*(\cO(3)\otimes \cL^{-1})$ has class $D_{n-1}$).
In the first case, there is a basis of the global sections of the line bundle defined by $D_{n+1}$ on $S$ such that the map
$$
\phi_{n+1}\phi_n^{-1}:\,S_n\,\longrightarrow\,S_{n+1}
$$
is given by any of the columns of the cofactor matrix:
$$
x\,\longmapsto (p_{n,1j}(x):\ldots:p_{n,4j}(x))\qquad(j\,=\,1,\ldots,4).
$$

\ts
In the proof of Proposition \ref{maps} we found that 
the divisor defined by $p_{n,ij}=0$ on $S_n$ is the sum of two effective divisors 
$\phi(C_{n,j})$, $\phi_n(C'_{n,i})$, 
where $C_{n,j}$ is the common zero locus of the $\phi_n^*p_{n,kj}$ 
and $C'_{n,i}$ is the common zero locus of the $\phi_n^*p_{n,ik}$ for $k=1,\ldots,4$.
These divisors both have genus $3$, so $(C_{n,j})^2=( C'_{n,i})^2=4$
and $C_{n,j}+C'_{n,i}=3D_n$.
%Their images under $\phi_n$ have degree $6$ in $\PP^3$, 
%so $D_n\cdot C_{n,j}=D_n\cdot C'_{n,i}=6$. 
As these divisors are effective, they are also ample (see Section \ref{effirr}), 
so their classes are in $\cA$.

%To find these classes, we observe that if $x\in \cA$ and $b(x,x)=4xx'=4$ 
%then $xx'=1$, so $x$ is a unit in the ring $\ZZ[\eta]$. 
%The group of units of $\ZZ[\eta]$ is  the set $\{\pm \eta^m:m\in\ZZ\}$. 
%As $b(\eta,\eta)=4\eta\eta'=-4$, 
%we have $x=\pm\eta^m$ with $m$ even  and the sign is $+$ since $x\in\cA$.
%Thus the two curves have class $D_k$, $D_{k'}$ for some $k,k'$ and %$D_k+D_{k'}=3D_n$.

Now we use that $D_n=\eta^{2n}$ and that multiplication by 
$\eta^{-2n}$ is an isometry of $N$ which maps $\cA$ into itself. 
%Thus we need to find $l\,(=k-n)$, $l'\,(=k'-n)$ such that $D_l+D_{l'}=3D_0$ and 
So we need to find $D:=a+b\eta$, $D':=a'+b'\eta\in\cA$ with sum $\eta^{-2n}(3D_n)=3D_0$. 
Hence $a+a'=3$, $b+b'=0$ and $a,a'>0$. 
Thus we may assume that $a=1$, $a'=2$. 
As $D_{l'}^2=4$, we then have $b'\neq 0$, hence also $b=-b'\neq 0$, 
and as $0<a^2+ab-b^2=1+b-b^2$, we get $b=1$. 
Hence $D=D_1$, $D'=D_{-1}$, so (up to permutation)
the class of $C_{n,j}$ is $\eta^{2n}D'=D_{n-1}$ and the class of
$C'_{n,i}$ is $\eta^{2n}D=D_{n+1}$. 
%(The identity
%$$
%\eta^4\,-\,3\eta^2\,+\,1\,=\,0\quad\mbox{implies}\quad
%\eta^{2n+2}\,-\,3\eta^{2n}\,+\,\eta^{2n-2}\,=\,0,\quad\mbox{hence}\quad
%3D_n\,=\,D_{n-1}\,+\,D_{n+1}~,
%$$
%which also shows that the sum of the two classes is $3D_n$.)

%(We observe that in Section \ref{quas} we showed that all divisors in any $|D_n|$
%are irreducible and reduced, hence $C_{n,i}$ and $C'_{n,j}$ are curves.)
As $C_{n,j}$ is the zero locus of a section of $\phi^*\cL$, the class of this line bundle is $D_{n-1}$. 
Therefore Proposition \ref{maps}, applied to ${}^tM_n$, 
shows that the columns of the cofactor matrix of $M_n$ give the map defined by the
global sections of the line bundle defined by $D_{n+1}$.
\qed

\

\section{Cayley's description of the automorphism}\label{cayleyA}

\subsection{} As observed by Cayley in \cite{C}, 
given a quartic determinantal surface, 
it is easy to find two others and to find isomorphisms between them. 
In our setup, starting from the projective model $S_n$ of $S$, 
he produces $S_{n+1}$ and $S_{n-1}$. 
The maps are those from Corollary \ref{class}.
Moreover, he shows that from the matrix $M_n$, whose determinant is the defining 
equation for $S_n$, one can find the matrices $M_{n-1},M_{n+1}$ which define 
$S_{n-1},S_{n+1}$ respectively.
The composition of the isomorphisms 
$S_0\rightarrow S_1\rightarrow S_2\rightarrow S_3=S_0$ is basically the automorphism $g$ we wanted to describe.

\subsection{A tritensor}\label{tritensor}
Let $S$ be a K3 surface with N\'eron Severi group 
$\NS(S)\cong (N,b)$ as in Section \ref{general}. 
Let $S_0=\phi_0(S)\subset \PP^3$, it is a smooth quartic determinantal surface. 
Let $M_0(x)$ be a $4\times 4$ matrix whose determinant defines $S_0$,
$$
S_0:\quad \det M_0(x)\,=\,0\qquad(\subset\PP^3)~.
$$
Write this matrix as
$$
M_0(x)\,:=\,\big(m_{kj}(x)\big)_{k,j=0,\ldots,3}~,
\qquad \mbox{with}\quad
m_{kj}(x)\,:=\,\sum_{i=0}^3 a_{ijk}x_i~,
$$
with coefficients $a_{ijk}\in\CC$.
The $4^3=64$ complex numbers $a_{ijk}$ can be viewed as the components of
a `tritensor' in $(\CC^4)^{\otimes 3}$. There are three obvious ways, up to transposition, in which this tritensor defines a $4\times 4$ matrix of linear forms. 
They are
$$
M_0(x)\,:=\,\big(m_{kj}(x)\big)~,\quad 
M_1(y)\,:=\,\big(m'_{ik}(y)\big)~,\quad
M_2(z)\,:=\,\big(m''_{ji}(z)\big)~,
$$
with coefficients 
$$
m_{kj}(x)\,:=\,\sum_{i=0}^3 a_{ijk}x_i~,\qquad
m'_{ik}(y)\,:=\,\sum_{j=0}^3a_{ijk}y_j~,\qquad
m''_{ji}(z)\,:=\,\sum_{k=0}^3a_{ijk}z_k~.
$$
The surprising thing is that the determinants of the $M_i$
define the quartic surfaces $S_i=\phi_i(S)$ respectively, 
provided these determinants are not identically zero, see below.
Thus the tritensor, or equivalently, any one of the matrices $M_i$, 
determines all the others. 
The maps between the $S_i$ are given by the rows and columns of the cofactor 
matrices of these matrices 
and they are thus defined by the tritensor as well.

\subsection{From $S_0$ to $S_{1}$} \label{cayley01}
Let $P_0(x)$ be the cofactor matrix of $M_0(x)$, so
$$
P_0(x) M_0(x)\,=\,M_0(x)P_0(x)\,=\,(\det M_0(x))I~,
$$
where $I$ is the $4\times 4$ identity matrix. 
After replacing $M_0$ and $P_0$ by their transposes if necessary, 
and after choosing a suitable basis of the global sections of the line bundle defined by $D_1$, 
we may assume (see Corollary \ref{class}) that the map
$\phi_1\phi_0^{-1}:S_0\rightarrow S_1$ is given by the columns of the adjoint matrix 
$P_0(x)$. 

On the other hand, for $x\in S_0$, we have $\det M_0(x)=0$ 
and thus each column of $P_0(x)$, 
provided it is not identically zero, provides a non-trivial solution to the linear
equations $M_0(x)y=0$. As $S_0$ is smooth, the rank of $M_0(x)$ is equal to three for any $x\in S_0$ and thus $y=y(x)$ is unique up to scalar multiple.

As $\det M_0(x)=0$ exactly for $x\in S_0$, we get, with some abuse of notation,
$$
S_1\,=\,\{y\,\in\,\PP^3:\;\exists x\in\,\PP^3\quad\mbox{s.t.}\quad   M_0(x)y\,=\,0\,\}~.
$$

The system of linear equations $M_0(x)y=0$ can be rewritten as:
$$
0\,=\,\sum_{j=0}^3\big(\sum_{i=0}^3 a_{ijk}x_i\big)y_j\,=\,
\sum_{i,j=0}^3 a_{ijk}x_iy_j\,=\,
\sum_{i=0}^3\,x_i\big(\sum_{j=0}^3 a_{ijk}y_j\big)\qquad(k=0,\ldots,3)~.
$$
This set of four equations is equivalent to the matrix equation, 
$$
{}^txM_1(y)\,=\,0,\qquad M_1(y)\,:=\,\big(m'_{ik}(y)\big),\quad m'_{ik}(y)\,:=\,\sum_{j=0}^3a_{ijk}y_j~.
$$
For $y\in\PP^3$ these equations have a non-trivial solution $x=x(y)$ if and only if
$\det M_1(y)=0$. In case $\det M_1(y)$ is not identically zero, 
it is a quartic polynomial that vanishes on the quartic surface $S_1$, 
and thus it is a defining equation for $S_1$. 
We will \emph{assume} that $ \det M_1(y)$ is not identically zero. Then
$$
S_1:\quad \det M_1(y)\,=\,0\qquad(\subset\PP^3)~.
$$
%and as $y=y(x)$ is a row of the cofactor matrix $P_0(x)$, we indeed have %$S_1=\phi_1\phi_0^{-1}(S_0)=\phi_1(S)$.

\subsection{From $S_1$ to $S_{2}$ and back to $S_0$} \label{cayley120}
We can repeat the procedure from section \ref{cayley01}: 
let $P_1(y)$ be the cofactor matrix of $M_1(y)$. 
As $S_0$ consists of the points $x$ with ${}^txM_1(y)=0$ and 
as $P_1(y)M_1(y)=0$ for $y\in S_1$, 
each row of $P_1$  defines the map $S_1\rightarrow S_0$. 
Thus the columns of $P_1$ define the map 
$\phi_2\phi_1^{-1}:S_1\rightarrow S_2=\phi_2(S)$
for a suitable basis of the global sections of the line bundle defined by $D_2$.
For $y\in S_1$, 
each column of $P_1(y)$ is then both a solution $z$ of $P_1(y)z=0$
and it defines a point of $S_2$. Hence we have
$$
S_2\,=\,\{z\,\in\,\PP^3:\;\exists y\in\,\PP^3\quad\mbox{s.t.}\quad   \;M_1(y)z\,=\,0\,\}~.
$$
Rewriting the linear equations $M_1(y)z=0$, we get ${}^tyM_2(z)=0$ with
$$
M_2(z)\,:=\,\big(m''_{ji}(z)\big),\quad m''_{ji}(z)\,:=\,\sum_{k=0}^3a_{ijk}z_k~.
$$

Now, \emph{assuming moreover} that $\det M_2(z)$ is not identically zero, 
we get
$$
S_2:\quad \det M_2(z)\,=\,0\qquad(\subset\PP^3)~.
$$

Finally we consider the cofactor matrix $P_2(z)$ of $M_2(z)$. 
The columns of $P_2(z)$ provide us with the map 
$\phi_3\phi_2^{-1}:S_2\rightarrow S_3=\phi_3(S)$,
for a suitable basis of the global sections of the line bundle defined by $D_3$.
Each column is also a solution of $M_2(z)x=0$. Rewriting this system, we get ${}^tzM_0(x)=0$, showing that $S_3$ is defined by $\det M_0(x)=0$ since this determinant \emph{is} not identically zero, being the defining equation of $S_0$. 

Thus $S_0=S_3$(!) and
the composition of the maps, each given by cubic polynomials (the minors of the matrices $M_i$)
$$
S_0\,\longrightarrow\,S_1\,\longrightarrow\,S_2\,\longrightarrow\,S_3\,=\,S_0
$$
is the map
$$
(\phi_3\phi_2^{-1})(\phi_2\phi_1^{-1})(\phi_1\phi_0^{-1})\,=\,\phi_3\phi_0^{-1}\,=\,
\phi_0 g\phi_0^{-1},
$$
where we used that $g=\phi_0^{-1}\phi_3$ (see Section \ref{autg})
and $g$ is the automorphism constructed by Oguiso.
To quote Cayley: ``The process may be indefinitely repeated".

\

\section{An explicit example}\label{explicit}

\subsection{The method}\label{exmet}
In this section we give an explicit example of a determinantal
K3 surface with Picard  number two. 
The main problem is to give an upper bound for the Picard number. For this we use 
a method described in \cite{RVL}. For the definition of the \'etale cohomology groups $\hH^i_{\et}(X, \QQ_\ell)$ and 
$\hH^i_{\et}(X, \QQ_\ell (1))$ for a scheme $X$, with values in $\QQ_\ell$ or its Tate twist $\QQ_\ell(1)$, we refer to \cite{T} and \cite{Mi}, p.\ 163--165.

The following result shows that if a smooth projective surface $X$ over a number field $K$ has good reduction at a prime $\p$, then the geometric Picard number of $X$ is bounded from above by the geometric Picard number of the reduction. 

\subsection{Proposition} \label{Inj}
Let $K$ be a number field with ring of integers $\cO$, let $\p$ be a prime of $\cO$ with residue field $k$, and let $\cO_\p$ be the localization of $\cO$ at $\p$. Let $\X$ be a smooth projective surface over $\cO_\p$ and set 
$X_{\Kbar} = \X \times_{\cO_\p} \Kbar$ and 
$X_{\kbar} = \X \times_{\cO_\p} \kbar$. 
Let $\ell$ be a prime not dividing $q=\# k$. 
Let $F_q^*$ denote the automorphism of 
$\hH^2_{\et}(X_{\kbar}, \QQ_\ell (1))$ induced by the $q$-th power Frobenius $F_q \in \Gal(\kbar/k)$.\\
Then there are natural injections
$$
\NS(X_{\Kbar})\otimes_{\ZZ} \QQ_\ell\, \inj\,
\NS(X_{\kbar})\otimes_{\ZZ} \QQ_\ell \,\inj\,
\hH^2_{\et}(X_{\kbar}, \QQ_\ell (1))~,
$$
that respect the intersection pairing and the action of Frobenius, respectively. 
The rank of $\NS(X_{\overline{k}})$ is at most the 
number of eigenvalues of $F_q^*$ that are roots of unity, counted with multiplicity.

\ts
The N\'eron-Severi group modulo torsion is isomorphic to the group of divisor classes modulo 
numerical equivalence (see \cite{F}, 19.3.1.(ii)). Therefore, the first injection, as well 
as the fact that it respects the intersection pairing, follows from 
\cite{F}, Examples 20.3.5 and 20.3.6. The second injection is in \cite{Mi}, Remark V.3.29.(d).
Each class $c \in \NS(X_{\kbar})$ is represented by 
a divisor, which is defined over some finite field; hence, some power of Frobenius 
fixes $c$. Since the N\'eron-Severi group $\NS(X_{\kbar})$ is finitely 
generated (see \cite{F}, 19.3.1.(iii)), it follows that some power of Frobenius acts
as the identity on $\NS(X_{\kbar})$. This implies the last statement. 
See also Proposition 6.2 and Corollary 6.4 in~\cite{RVL1} (which counts the eigenvalues 
that are a root of unity up to a factor $q$ because it refers to the action on 
$\hH^2_{\et}(X_{\kbar}, \QQ_\ell)$ without the Tate twist).
\qed

\subsection{Proposition} \label{Lefschetz}
Let $X$ be a K3 surface over a finite field $k\cong \FF_q$. 
As in Proposition \ref{Inj}, let $F_q^*$ denote the automorphism of 
$\hH^2_{\et}(X_{\kbar}, \QQ_\ell (1))$ induced by the $q$-th power Frobenius $F_q \in \Gal(\kbar/k)$, 
and for any $n$, let $\Tr((F_q^*)^n)$ denote the trace of $(F_q^*)^n$. Then we have 
$$
\Tr\big((F_q^*)^n\big) = \frac{\#X(\FF_{q^n})-1-q^{2n}}{q^n}.
$$
Furthermore, the characteristic polynomial $f(t) = \det (t-F_q^*) \in \QQ[t]$ of $F_q^*$ has degree $22$ and satisfies the functional equation 
$$
t^{22} f(t^{-1}) = \pm f(t).
$$

\ts 
Let $F_X$ be the $q$-th power absolute Frobenius of $X$, which acts as the identity on points and by raising to the $q$-th power on the coordinate rings of affine open subsets of $X$. The geometric Frobenius $\varphi = F_X \times 1$ on $X \times_k \kbar = X_{\kbar}$ is an endomorphism of $X_{\kbar}$ over $\kbar$ (cf.\ \cite{Mi}, proof of V.2.6 and pages 290--291). The set of fixed points of $\varphi^n$ is $X(\FF_{q^n})$. The Weil conjectures (see \cite{Mi}, \S VI.12) state that the eigenvalues of $\varphi^*$ acting on $\hH^i_{\et}(X_{\kbar}, \QQ_\ell)$ have absolute value $q^{i/2}$. Since $X$ is a K3 surface, we have $\dim \hH^i_{\et}(X_{\kbar}, \QQ_\ell)=1,0,22,0,1$ for $i=0,1,2,3,4$, respectively (see \cite{badescu}, 8.4 and Theorem 10.3), so the Lefschetz trace formula for $\varphi^n$ (see \cite{Mi}, Theorems VI.12.3 and VI.12.4) yields
\begin{equation}\label{simplifiedlefschetz}
\# X(\FF_{q^n}) = \sum_{i=0}^4 (-1)^i \Tr\big((\varphi^*)^n|\hH^i_{\et}(X_{\kbar},\QQ_\ell)\big)= 1 + \Tr\big((\varphi^*)^n|\hH^2_{\et}(X_{\kbar},\QQ_\ell)\big) + q^{2n}.
\end{equation}
For the remainder of this proof we restrict our attention to the middle cohomology, so $\hH^i_{\et}$ with $i=2$. 
By the Weil conjectures, the characteristic polynomial $f_\varphi(t) = \det (t-\varphi^*|\hH^2_{\et}(X_{\kbar}, \QQ_\ell))$ 
is a polynomial in $\ZZ[t]$ satisfying the functional equation $t^{22}f_\varphi(q^2/t) = \pm q^{22} f_\varphi(t)$ (note that the polynomial
$P_2(X,t) = \det(1-\varphi^* t|\hH^2_{\et}(X_{\kbar}, \QQ_\ell))$ of \cite{Mi}, \S VI.12, is the reverse of $f_\varphi$). 
Let $\varphi^*(1)$ denote the action on $\hH^2_{\et}(X_{\kbar},\QQ_\ell(1))$ (with a Tate twist) induced by $\varphi$.
Note that the fact that $\varphi^*(1)$ acts on the middle cohomology is not reflected in the notation. 
The eigenvalues of $\varphi^*(1)$ differ from those of $\varphi^*$ on $\hH^2_{\et}(X_{\kbar},\QQ_\ell)$ by a factor $q$ (see \cite{T}), so we have 
\begin{equation}\label{factorq}
\Tr\big((\varphi^*)^n|\hH^2_{\et}(X_{\kbar},\QQ_\ell)\big) = q \cdot \Tr\big(\varphi^*(1)^n\big),
\end{equation}
and the characteristic polynomial $f_\varphi^{(1)}\in \QQ[t]$ of $\varphi^*(1)$ satisfies 
$q^{22}f_\varphi^{(1)}(t) = f_\varphi(qt)$, and thus the functional equation $t^{22} f_\varphi^{(1)}(1/t) = \pm f_\varphi^{(1)}(t)$.
It follows that the eigenvalues, and hence the characteristic polynomials, of $\varphi^*(1)$ and $\varphi^*(1)^{-1}$ coincide. 
Finally, the product of the geometric Frobenius $\varphi = F_X \times 1$ and the Galois automorphism $1 \times F_q$ on $X\times_k \kbar = X_{\kbar}$ is the absolute Frobenius $F_{X_{\kbar}}$, which acts as the identity on the cohomology groups, so the maps $\varphi^*(1)$ and $F_q^*$ act as inverses of each other (see \cite{Mi}, Lemma VI.13.2 and Remark VI.13.5, and \cite{T}, \S 3). 
We conclude $f = f_\varphi^{(1)}$ and $\Tr\big((F_q^*)^n\big) = \Tr (\varphi^*(1)^{-n}) = \Tr (\varphi^*(1)^{n})$,
which, together with \eqref{simplifiedlefschetz} and \eqref{factorq}, implies the proposition.
\qed

\subsection{An explicit determinant}\label{exdet}
Consider the following matrix, whose entries are linear forms in the variables 
$x_0,\ldots,x_3$ with integer coefficients:
%\begin{equation}\label{eq:defN}
$$
M_0=
\begin{pmatrix}
x_0 & x_2 & x_1+x_2 & x_2+x_3 \\
x_1 & x_2+x_3 & x_0+x_1+x_2+x_3 & x_0+x_3\\
x_0+x_2 & x_0+x_1+x_2+x_3 & x_0+x_1 & x_2\\
x_0+x_1+x_3 & x_0+x_2 & x_3 & x_2
\end{pmatrix}~.
$$
%\end{equation}

The following theorem shows that any such matrix that is congruent to $M_0$ modulo $2$ 
defines a determinantal quartic surface in $\PP^3=\PP^3(\CC)$ with Picard number $2$.

\subsection{Theorem} \label{PicNum}
Let  $R=\ZZ [ x_0, x_1, x_2, x_3]$ and let $M \in M_4(R)$ be any matrix whose entries are homogeneous polynomials of degree $1$ 
and such that $M$ is congruent modulo $2$ to the matrix $M_0$ above. 
Denote by $S$ the complex surface in $\PP^3$ given by $\det M=0$. 
Then $S$ is a K3 surface and its Picard number equals 2.

\ts
Let $\fS$ denote the surface over the localization $\ZZ_{(2)}$ of $\ZZ$ at the prime $2$ given by $\det M=0$, and 
write $S'$ and $\overline{S'}$ for $\fS_{\FF_2}$ and 
$\fS_{\overline{\FF}_2}$, respectively. 
One checks that $S'$ is smooth and $\fS$ is reduced, for instance with Magma \cite{magma}. 
Since $\Spec \ZZ_{(2)}$ is integral and regular of dimension $1$, the scheme $\fS$ is integral, and 
the map $\fS \to \Spec \ZZ_{(2)}$ is dominant, it follows from \cite{H}, Proposition III.9.7, that $\fS$ is flat over 
$\Spec \ZZ_{(2)}$. Since the fiber over the closed point is smooth, it follows from \cite{L}, Definition 4.3.35, that $\fS$
is smooth over $\Spec \ZZ_{(2)}$. Therefore, $S=\fS_\CC$ is smooth as well, so $S$ and $S'$ are K3 surfaces. 
Let $F_2^*$ denote the automorphism of 
$\hH^2_{\et}(\overline{S'}, \QQ_\ell (1))$ induced by Frobenius $F_2 \in \Gal(\overline{\FF}_2/\FF_2)$.

The divisor classes in $\hH^2_{\et}(\overline{S'}, \QQ_\ell (1))$ 
defined by the hyperplane class and the curve $C$ as in Proposition \ref{line} 
span a two-dimensional subspace $V$ 
on which $F_2^*$ acts as the identity.
We denote the linear map induced by $F_2^*$ on the quotient $W:=\hH^2_{\et}(\overline{S}_2, \QQ_\ell(1) )/V$
by $F_2^*|_W$, so that $\Tr (F_2^*)^n = \Tr (F_2^*|_{V})^n+ \Tr (F_2^*|_W)^n = 2+\Tr (F_2^*|_W)^n$ for every integer $n$.
From Proposition \ref{Lefschetz}, we obtain
$$
\Tr (F_2^*|_W)^{n} \,=\, -2\,+\,\frac{\# S'(\FF_{2^n})-1-2^{2n}}{2^n}~.
$$
We counted the number of points in $S'(\FF_{2^n})$ for $n=1,\ldots ,10$ with Magma \cite{magma}. 
The results are in the table below.
{\renewcommand{\arraystretch}{1.4}
$$
\begin{array}{|c|c|c|c|c|c|c|c|c|c|c|}\hline
n&1&2&3&4&5&6&7&8&9&10 \\
\hline
\# S'(\FF_{2^n}) &6&26&90&258&1146&4178&17002&64962&260442&1044786 \\ \hline
\Tr(F_2^*|_W)^n &-\frac{3}{2} & \frac{1}{4} & \frac{9}{8} & -\frac{31}{16} & \frac{57}{32} & -\frac{47}{64} & \frac{361}{128} & -\frac{1087}{256} & -\frac{2727}{512} & -\frac{5839}{1024} \\
\hline
\end{array}
$$
}
If $\lambda_1,\ldots , \lambda_{20}$ denote the eigenvalues of $F_2^*|_W$, 
then the trace of $(F_2^*|_W)^n$ equals 
$$
\mathrm{Tr}(F_2^*|_W)^{n} \,=\,\lambda_1^n+\ldots + \lambda_{20}^n~,
$$
i.e., it is the $n$-th power sum symmetric polynomial in the eigenvalues of $F_2^*|_W$. 
Let $e_n$ denote the elementary symmetric polynomial of degree $n$ in the 
eigenvalues of $F_2^*|_W$ for $n \geq 0$. Using Newton's identities
$$
ne_n = \sum_{i=1}^n (-1)^{i-1} e_{n-i} \cdot \mathrm{Tr}(F_2^*|_W)^{i}
$$
and $e_0=1$, we compute the values of $e_n$ for $n=1,\ldots ,10$. They are 
listed in the following table.

{\renewcommand{\arraystretch}{1.4}
$$
\begin{array}{|c|c|c|c|c|c|c|c|c|c|c|}\hline
n&1&2&3&4&5&6&7&8&9&10 \\
\hline
e_n&-\frac{3}{2}&1&0&0&0&0&\frac{1}{2}&0&-1&2\\
\hline
\end{array}
$$
}

We denote the characteristic polynomial of a linear operator $T$ by $f_T$, so that 
$$
f_{F_2^*} = f_{F_2^*|_{V}}\cdot f_{F_2^*|_W}  = (t-1)^2 f_{F_2^*|_W}~.
$$
Because $f_{F_2^*}$ satisfies the functional equation of Proposition \ref{Lefschetz},
the polynomial $f_{F_2^*|_W}$ satisfies $t^{20}f_{F_2^*|_W}(t^{-1}) =  \pm f_{F_2^*|_W}(t)$.
Since the middle coefficient $e_{10}=2$ of $t^{10}$ in $f_{F_2^*|_W}$ is nonzero, the sign in this 
functional equation is $+1$, so $f_{F_2^*|_W}$ is palindromic and we get 
\begin{align*}\label{fF2W}
f_{F_2^*|_W}  &=t^{20}-e_1t^{19} +e_2 t^{18} - \dots +e_{10}t^{10} -e_9t^9 + \dots -e_1t + 1  \\
& = t^{20}+\tfrac{3}{2}t^{19}+t^{18}-\tfrac{1}{2}t^{13}+t^{11}+2t^{10}+t^9-\tfrac{1}{2}t^7+t^2+\tfrac{3}{2}t+1.
\end{align*}
With Magma \cite{magma}, one checks that this polynomial is irreducible over $\QQ$, and as it is not integral, its roots 
are not algebraic integers, so none of its roots is a root of unity.
Hence, the polynomial $f_{F_2^*} = (t-1)^2f_{F_2^*|_W}$ has exactly two roots that are a root of unity. 
This implies that $F_2^*$  has only two eigenvalues (counted with multiplicity) 
that are roots of unity, and so, by Proposition~\ref{Inj}, it follows 
that the rank of the Picard group $\NS(S) \cong \NS(\fS_{\bar{\QQ}})$ is bounded by two from above.
On the other hand, by Proposition \ref{propB} we know that the rank is at least two,
hence $\NS(S)$ has rank two.
\qed

\plotpicture

\subsection{The automorphism}\label{exaut}
We apply Cayley's method, starting with the choice $M_0(x)$ as in Section \ref{exdet}. 
The matrices $M_1(y)$ and $M_2(z)$ as in Section \ref{tritensor} are given by 
$$
M_1(y)=
\begin{pmatrix}
 y_0&y_2 + y_3&y_0 + y_1 + y_2& y_0 + y_1\\
 y_2&y_0 + y_2&y_1 + y_2&y_0\\
y_1 + y_2 + y_3&y_1 + y_2&y_0 + y_1 + y_3&y_1 + y_3\\
y_3&y_1 + y_2 + y_3&y_1&y_0 + y_2
\end{pmatrix}~,
$$
and 
$$
M_2(z)=
\begin{pmatrix}
z_0 + z_2 + z_3&z_1 + z_3&z_2& z_3\\
z_2 + z_3&z_2&z_0 + z_1 + z_2 + z_3&z_1 + z_2\\
z_1 + z_2&z_0 + z_1 + z_2&z_0 + z_1&z_1 + z_3\\
z_1& 0& z_0 + z_2 + z_3& z_0 + z_1
\end{pmatrix}~.
$$
The determinants of $M_1(y)$ and $M_2(z)$ are not identically zero, so for any $n \in \ZZ$, the surface
$S_n$ is given by the vanishing of $\det M_{(n \!\! \mod 3)}$. 
As described in Sections \ref{cayley01} and \ref{cayley120}, the map $S_n \to S_{n+1}$ is given by any column 
of the cofactor matrix $P_{n'}$ of $M_{n'}$ with $n' = n \!\! \mod 3$. For example, if we write 
$M_0(x) = \big(m_{kj}(x)\big)_{k,j=0,\ldots,3}$ as before, then the first column of the cofactor matrix $P_0(x)$
is $(g_0,g_1,g_2,g_3)^{\rm t}$ with 
\begin{align*}
  g_0 =&  -x_0^3 - x_0^2x_1 - x_0^2x_2 - 2x_0x_1x_2 + x_0x_2^2 + 2x_0x_3^2 - x_1^2x_2 - 2x_1x_2x_3 + x_1x_3^2 \\
       &   - 2x_2^2x_3 - x_2x_3^2 + x_3^3, \\
  g_1 =&  x_0^3 + 2x_0^2x_1 + x_0^2x_3 + x_0x_1^2 - 2x_0x_1x_2 + 3x_0x_1x_3 + x_0x_2^2 - 2x_0x_2x_3 - 2x_1^2x_2 \\
       & + x_1^2x_3 - x_1x_2x_3 + x_1x_3^2 + x_2^3 - 2x_2x_3^2, \\
  g_2 =&  -2x_0^2x_1 + x_0^2x_2 - 2x_0^2x_3 - x_0x_1^2 - x_0x_1x_2 - 4x_0x_1x_3 + x_0x_2^2 - 3x_0x_3^2 + x_1^2x_2\\
       & - x_1^2x_3 + x_1x_2^2 + x_1x_2x_3 - 2x_1x_3^2 - x_2^3 + x_2^2x_3 - x_3^3, \\
  g_3 =&  3x_0^2x_1 - 2x_0^2x_2 + x_0^2x_3 + 4x_0x_1^2 + x_0x_1x_2 + 3x_0x_1x_3 - 2x_0x_2^2 + 2x_0x_2x_3 + 3x_0x_3^2 \\
       & + x_1^3 + 2x_1^2x_2 + x_1^2x_3 + 2x_1x_2x_3 + x_1x_3^2 - x_2^3 + x_2^2x_3 + 3x_2x_3^2 + x_3^3.
\end{align*}
Hence, the map $S_0 \to S_1$ is given, at least on an open subset of $S_0$, by sending 
$[x_0:x_1:x_2:x_3]$ to $[g_0(x):g_1(x):g_2(x):g_3(x)]$.
However, the polynomials $g_0,g_1,g_2,g_4$ all vanish on a curve $C$ as described in Corollary \ref{class}. In 
order to define the map $S_0 \to S_1$ everywhere, we use the other three columns of the cofactor matrix $P_0(x)$. 
Similarly, the columns of the cofactor matrices $P_1(y)$ and $P_2(z)$ of $M_1(y)$ and $M_2(z)$ determine the maps 
$S_1 \to S_2$ and $S_2 \to S_0$, respectively. 

Cayley's method therefore gives the automorphism $g$ explicitly as the composition $S_0 \to S_1 \to S_2 \to S_0$ 
of three maps, each given by cubic polynomials. We can thus describe $g$ explicitly by quadruples of coordinate functions, 
each of which homogeneous of degree $3^3=27$ in $x_0,x_1,x_2,x_3$. Unfortunately, these quadruples are far too large 
to write down.

\subsection{Defining polynomials of lower degree} \label{deg18}
We showed in section \ref{findg} that the automorphism $g$ can in fact be given by polynomials $R_{ij}$ of degree $18$. 
We now describe how we used linear algebra to find such polynomials explicitly. Unfortunately, these polynomials are
still far too large to publish, but they are electronically available from the authors upon request. 

We are looking for quadruples $(G_0,G_1,G_2,G_3) \in \QQ[x_0,x_1,x_2,x_3]^4$ of homogeneous polynomials of degree $18$, 
such that the rational map $\PP^3 \dashrightarrow \PP^3$ given by $x \mapsto [G_0(x):G_1(x):G_2(x):G_3(x)]$ coincides 
on an open subset of $S=S_0$ with our automorphism $g$. There are 
$\binom{18+3}{3}=1330$ monomials of degree $18$, 
so this gives $4\cdot 1330 = 5320$ unknown coefficients. For each point $x = [x_0:x_1:x_2:x_3] \in S$, we can compute
$g(x) = [y_0:y_1:y_2:y_3]$ with Cayley's method above. The identity $g(x) = [G_0(x): G_1(x): G_2(x): G_3(x)]$ is equivalent 
with the six equalities $y_jG_i(x) = y_iG_j(x)$ for $0\leq i < j \leq 3$, which are linear in the unknown coefficients of 
$G_0, G_1,G_2,G_3$. In fact, if $x$ is defined over a field extension $K$ of $\QQ$ of degree $d$, then we can write 
the equation $y_jG_i(x) = y_iG_j(x)$ in terms of a basis for $K$ over $\QQ$; since we know that the $G_i$ can be chosen 
over $\QQ$, we can split up the equation in $d$ independent equations, each linear in the unknown coefficients. 
After choosing sufficiently many points over various number fields, we obtain a large system of equations over $\QQ$ that we solved with Magma \cite{magma}. The solution space $V$ has dimension $2724$ inside the space of quadruples 
$(G_0,G_1,G_2,G_3) \in \QQ[x_0,x_1,x_2,x_3]^4$ of homogeneous polynomials of degree $18$. This space contains the 
space $U$ of quadruples that vanish on $S$, which has dimension $4 \cdot \binom{17}{3} = 2720$. 

To verify that we used enough points, we took four quadruples $(G_{s1},G_{s2}, G_{s3},G_{s4})$, for $s=1,2,3,4$, that generate the quotient $V/U$ and checked that they indeed define the same map 
as $g$ on some open subset of $S$; this can be done by taking, as in the previous section, a quadruple 
$(F_0,F_1,F_2,F_3)$ of homogeneous polynomials of degree $27$ describing $g$, and checking that for each $1 \leq s \leq 4$, and each $0\leq i<j \leq 3$, the polynomial $F_iG_{sj}-F_jG_{si}$ is divisible 
by the defining polynomial $\det M_0(x)$ of $S=S_0$. 

We also verified with Magma \cite{magma} that one can in fact choose three quadruples such that the $3\cdot 4$ polynomials in them have no common base points on $S$, i.e., 
at each point on $S$, the automorphism $g$ is defined by at least one of these three quadruples. This computation 
was done over the rational numbers and therefore holds over any field of characteristic $0$. 

\subsection{Points of period two} 

As mentioned in Section \ref{topLnum}, the topological Lefschetz number of $g^2$ is $344$. 
Using the previous section, we can explicitly define a scheme $\Xi$ over $\ZZ$ such that 
$\Xi_\QQ$ is $0$-dimensional and consists of points of period $2$ as follows. First we choose three quadruples 
$(G_{s1},G_{s2}, G_{s3},G_{s4}) \in \ZZ[x_0,x_1,x_2,x_3]^4$, for $s=1,2,3$, of homogeneous polynomials of 
degree $18$ that together describe $g$ everywhere on $S$. Similarly, we compute three quadruples 
$(H_{s1},H_{s2}, H_{s3},H_{s4}) \in \ZZ[x_0,x_1,x_2,x_3]^4$, for $s=1,2,3$, of homogeneous polynomials of 
degree $18$ that together describe the inverse 
$g^{-1}$ everywhere on $S$; as described in Section \ref{cayley120}, this is done by using the rows of the cofactor matrices
rather than the columns. A point $x \in S$ has period $2$ if and only if $g(x) = g^{-1}(x)$, so if and only if 
$H_{si}(x)G_{tj}(x) = H_{sj}(x)G_{ti}(x)$ for all $1\leq s,t \leq 3$ and all $0 \leq i < j \leq 3$. Hence, these $54$ equations 
of degree $36$, together with the defining polynomial $\det M_0$ for $S$, define a scheme $\Xi$ over $\ZZ$ such that the scheme $\Xi_\QQ$ consists of all points of period $2$. 

The following proposition states that in our explicit example, none of these points is defined over a number field 
of degree less than $344$ over $\QQ$. The proof is based on reduction modulo primes of good reduction, as we were unable to 
perform any significant computations with $\Xi$ in Magma \cite{magma} over the rational numbers. All computations 
are available from the authors upon request. 

\subsection{Proposition}\label{periodic}
Let $\Pi \subset S(\bar{\QQ})$ denote the set of all points of period $2$ under the automorphism $g$.
Then $\#\Pi = 344$ and the Galois group $\Gal(\bar{\QQ}/\QQ)$ acts transitively on $\Pi$.

\ts  Let $K$ be a finite Galois extension of $\QQ$ with $\Pi \subset S(K)$ and let $\cO_K$ be the ring of integers of $K$. Note that $\Pi$ consists of the $K$-points, or equivalently, the $\cO_K$-points, of the $0$-dimensional scheme $\Xi$ 
constructed above. Take $p \in \{17,101\}$, and let $\p \subset \cO_K$ be a prime ideal above prime $p$. 
Set $k(\p) = \cO_K/\p$ and let $\bar{\FF}_p$ be an algebraic closure of $k(\p)$.  
One checks with Magma \cite{magma} that the surface $S_{\ZZ} \subset \PP^3_{\ZZ}$ given by $\det M_0=0$ has good reduction modulo $p$, i.e., the reduction $S_p = S_\ZZ \times \FF_p$ is smooth, and the same holds for the surfaces in $\PP^3_{\ZZ}$ given by $\det M_i=0$ for $i=1,2$. We claim that the composition $\Pi= \Xi(K) = \Xi(\cO_K) \to \Xi(k(\p)) \to \Xi(\bar{\FF}_p)$ of the reduction map and the inclusion is surjective. Indeed, one can verify with Magma \cite{magma} 
that $\Xi_{\FF_p}$ has dimension $0$, degree $344$, and is reduced. Since $\FF_p$ is perfect, it follows that $\Xi_{\bar{\FF}_p}$ is reduced (see \cite{egaiv2}, Proposition 4.6.1). Let $\Xi_0$ be an irreducible component of $\Xi_{\bar{\FF}_p}$. Then $\Xi_0$ is $0$-dimensional, so it is affine and its coordinate coordinate ring $A(\Xi_0)$ is $0$-dimensional and Noetherian, and therefore Artin (see \cite{AM}, Theorem 8.5). Hence, $A(\Xi_0)$ is the product of local Artin rings (\cite{AM}, Theorem 8.7), and since it also integral, it is local itself. From \cite{AM}, Proposition 8.6, and the fact that $A(\Xi_0)$ is integral, we conclude that $A(\Xi_0)$ is a field, which, being a finite extension of $\bar{\FF}_p$, is isomorphic to $\bar{\FF}_p$. Thus, $\Xi_0$ is a smooth point, and $\Xi_{\FF_p}$ is smooth over $\FF_p$. Hence, it follows from Hensel's Lemma
that every point of $\Xi_{\FF_p}$ over some finite extension of $\FF_p$ lifts to some finite extension of $\ZZ_p$ and, since every point on a $0$-dimensional scheme is algebraic, to some finite extension of $\QQ$.  
As the topological Lefschetz number of $g^2$ equals $344$, we find that $g$ has at most $344$ points of period $2$, so 
$\# \Pi \leq 344$. From the claim and the equality $\# \Xi(\bar{\FF}_p) = \deg \Xi_{\FF_p} = 344$ we conclude that
$\# \Pi = 344$ and the reduction map $r\colon \Xi(K) \to \Xi(k(\p)) = \Xi(\bar{\FF}_p)$ is a bijection.

The bijection $r$ respects the Galois action of the decomposition subgroup $G_\p \subset \Gal(K/\QQ)$ associated to $\p$. 
Each Galois orbit $C$ of $\Xi(K)$ under $\Gal(K/\QQ)$ splits as the disjoint union of orbits under $G_\p$; since $G_\p$ naturally surjects onto $\Gal(k(\p)/\FF_p)$, the image $r(C) \subset \Xi(k(\p))$ splits as the disjoint union of orbits of 
$\Xi(k(\p))$ under $\Gal(k(\p)/\FF_p)$, or equivalently, of orbits of $\Xi(\bar{\FF}_p)$ under $\Gal(\bar{\FF}_p/\FF_p)$. This implies that 
the sizes of the Galois orbits of $\Xi(\bar{\FF}_p)$ form a partition of $344$ that is a refinement of the partition corresponding to the sizes of the Galois orbits of $\Xi(K)$ under $\Gal(K/\QQ)$. 
More precisely, if $m_1, \ldots, m_s$ are the sizes of the 
Galois orbits of $\Xi(K)$ under $\Gal(K/\QQ)$, and $n_1, \ldots, n_t$ are the sizes of the Galois orbits of $\Xi(\bar{\FF}_p)$ under $\Gal(\bar{\FF}_p/\FF_p)$, then $m_1+\dots + m_s = 344 = n_1+\dots +n_t$, and there is a partition $(I_1,\ldots,I_s)$ of the set $\{1, 2, \ldots, t\}$ such that $m_j = \sum_{i \in I_j} n_i$ for all $1 \leq j \leq s$. 

We computed the number of $\FF_{p^t}$-points on $\Xi_{\FF_p}$ for $p=17$ and $p=101$ and $ 1\leq t \leq 150$ with Magma \cite{magma}. 
For $p=17$ we found that there are $4$ points on $\Xi_{\FF_p}$ over $\FF_p$, as well as $2$ more points over $\FF_{p^2}$ and $12$ more points over $\FF_{p^{12}}$, which are not defined over a smaller field, $110$ more points over $\FF_{p^{55}}$, and no other points over any field $\FF_{p^t}$ with $t\leq 150$. It follows that $\Xi(\bar{\FF}_{17})$ has Galois orbits of sizes $1,1,1,1,2,12,55,55$, and as none of the remaining $216$ points is defined over a field of degree less than $150$ over $\FF_p$, one orbit of size $216$. 
For $p=101$ we found that there are $20$ points on $\Xi_{\FF_p}$ over $\FF_{p^{20}}$, which are not defined over a smaller field, $26$ more points over $\FF_{p^{26}}$, and no other points over any field $\FF_{p^t}$ with $t\leq 150$. It follows that $\Xi(\bar{\FF}_{101})$ has Galois orbits of size $20,26$, and as none of the remaining $298$ points is defined over a field of degree less than $150$ over $\FF_p$, one orbit of size $298$. 
The only partition of $344$ of which both the partitions $\{ 1,1,1,1,2,12,55,55,216 \}$ and $\{20,26,298\}$ are a refinement is the trivial partition $\{344\}$ of one part, so we find that $\Xi(K)$ is one orbit under $\Gal(K/\QQ)$, which proves the proposition. \qed

\end{document}